\documentclass[12pt]{article}
\usepackage[american]{babel}
\usepackage{amsmath}
\usepackage{latexsym}
\usepackage{amssymb}
\usepackage{theorem}
\usepackage{diagrams}
\theoremstyle{change}
\newtheorem{Thm}{Theorem}[section]
\newtheorem{Cor}[Thm]{Corollary}
\newtheorem{Prop}[Thm]{Proposition}
\newtheorem{Lem}[Thm]{Lemma}
{\theorembodyfont{\rmfamily}
\newtheorem{Num}[Thm]{}

}
\newcommand{\bra}[1]{\langle#1\rangle}
\newcommand{\proof}{\par\medskip\rm\emph{Proof. }}

\newcommand{\qed}{\ \hglue 0pt plus 1filll $\Box$}
\newcommand{\oot}{\longleftarrow}
\newcommand{\too}{\longrightarrow}
\newcommand{\mapstoo}{\longmapsto}

\newcommand{\RR}{\mathbb{R}}
\newcommand{\ZZ}{\mathbb{Z}}
\newcommand{\NN}{\mathbb{N}}

\newcommand{\QQ}{\mathbb{Q}}

\newcommand{\SKIP}[1]{}
\newcommand{\SO}{\mathrm{SO}}
\newcommand{\SL}{\mathrm{SL}}

\renewcommand{\emptyset}{\varnothing}
\newcommand{\tr}{\mathrm{tr}}
\newcommand{\supp}{\mathrm{supp}}

\begin{document}

\title{\bf Affine $\Lambda$-buildings,
ultrapowers of Lie groups and Riemannian symmetric spaces:
an algebraic proof of the Margulis conjecture}
\author{Linus Kramer and Katrin Tent
\thanks{Both supported by {\em Heisenberg} fellowships of the DFG}\\
{\small Universit\"at W\"urzburg,
Am Hubland,
D--97074 W\"urzburg,
Germany}\\
{\small {\tt \{kramer,tent\}@mathematik.uni-wuerzburg.de}}\\
\it Dedicated to the memory of Reinhold Baer on his 100th birthday}
\date{}
\maketitle

\section*{Introduction}
Let $X$ and $X'$ be Riemannian symmetric spaces of noncompact
type, and let $G=I(X)$ and $G'=I(X')$ denote their respective
isometry groups. Then
$G$ and $G'$ are (finite extensions of)
semisimple real Lie groups. Suppose furthermore that $\Gamma$
is a (torsion free) group which injects as a cocompact lattice into
$G$ and into $G'$,
\[
G\lInto\Gamma\rInto G'.
\]
Then $\Gamma$ acts cocompactly on $X$ and $X'$, and the orbit spaces
$X/\Gamma$ and $X'/\Gamma$ are closed manifolds with contractible
universal covers. The only
nontrivial homotopy groups are thus the fundamental groups,
\[
\pi_1(X/\Gamma)\cong\Gamma\cong\pi_1(X'/\Gamma).
\]
It follows that there is
a homotopy equivalence $X/\Gamma\rTo^\simeq X'/\Gamma$ which
can be lifted to a map $f:X\rTo X'$. This map $f$ can then be shown to be
a \emph{quasi-isometry}. Recall that a (not necessarily continuous)
map $f:X\too X'$ between metric spaces if called
an $(L,C)$-quasi-isometry if there exist constants $L\geq 1$ and $C\geq0$
such that
\[
L^{-1}\cdot d(x,y)-C\leq d'(f(x),f(y))\leq L\cdot d(x,y)+C.
\]
holds for all $x,y\in X$, and such that every $x'\in X'$ is $C$-close to
some image point $f(x)$. Special instances of quasi-isometries are

(1) \emph{isometries} (where $(L,C)=(1,0)$)

(2) \emph{bi-Lipschitz maps} (where $C=0$)

(3) \emph{coarse isometries} (where $L=1$).

\noindent

\smallskip
\noindent\textbf{Mostow's Rigidity Theorem} \cite{Mostow}
\emph{If the Riemannian symmetric spaces
$X$ and $X'$ have no de Rham factors of rank $1$, then there
exists a Lie group isomorphism $G\rTo^\cong G'$ such the diagram
\begin{diagram}[size=2em]
\Gamma & \rInto & G \\
\dInto & \ldTo^\cong  \\
G' 
\end{diagram}
commutes.
In particular, there is an isometry $X\cong X'$
(after rescaling the metrics on the de Rham factors of $X'$, if necessary.)}

\medskip
The key idea of Mostow's proof was to show that $f$ induces an
isomorphism between the spherical buildings attached to $X$ and
$X'$. This is not trivial, and the group $\Gamma$ plays an important
r\^ole in Mostow's
argument. The Margulis conjecture may be stated as follows:

\smallskip\noindent
\textbf{Margulis Conjecture}
\emph{If $f:X\too X'$ is a quasi-isometry, then (after rescaling the
de Rham factors) $X$ and $X'$ are isometric, and there exists
an isometry $\bar f:X\too X'$ at bounded distance from $f$.}

\medskip
Kleiner-Leeb \cite{KL}
proved the Margulis conjecture in 1997;
later, Eskin-Farb \cite{EF} gave another proof. Kleiner-Leeb used
Gromov's technique of asymptotic cones \cite{Gromov}.
The asymptotic cone is a very powerful functor on the category of 
pointed metric spaces. Kleiner-Leeb
showed that the asymptotic
cone of a Riemannian symmetric space is an affine $\RR$-building, and
they used this building systematically in their work \cite{KL}.
(Note that they work with a class of metric spaces they call
\emph{euclidean buildings}; their system of axioms is different
from Tits' \emph{syst\`emes d'appartements} \cite{Tits}.
Parreau \cite{Parreau} proved that their spaces are a subcase of
Tits' syst\`emes d'appartements.)
The point is that the quasi-isometry between the Riemannian
symmetric spaces
induces a homeomorphism of their respective affine $\RR$-buildings, and
that this homeomorphism can be shown to be an isomorphism.

However, it does not follow from their paper which building the asymptotic
cone of a given Riemannian symmetric space really is.
In this paper, we give a
general group-theoretic 
construction of affine $\RR$-buildings, and more generally, of
affine $\Lambda$-buildings, associated to semisimple Lie groups over
nonarchimedean real closed fields. This result is new.
The construction of Kleiner-Leeb
using the asymptotic cone appears only as a special case.
Also, we give a new, sheaf-theoretic proof for the topological
rigidity of affine $\RR$-buildings.
The explicit knowledge of the building arising here as the asymptotic
cone simplifies
the proof of the Margulis conjecture; it also sheds some light
on the algebraic background of Mostow rigidity and the Margulis conjecture.

The present report is meant as a survey without detailed proofs.
After completing it, we learned that Thornton, a student of Kleiner,
independently proved Corollary \ref{HomogeneousCone} \cite{Thornton}.
He writes in his introduction: 'In the case of [...] the group
$\SL_n\RR$, Leeb identified the asymptotic cone as a homogeneous
space over an algebraic group. Parreau showed
that the asymptotic cone for [...] $\SL_n\RR$ fits a certain model for
Euclidean buildings [...]'. The case $G=\SL_n\RR$ is in fact rather
elementary; Bennett \cite{Bennett} showed that one obtains
from this group an affine $\Lambda$-building over \emph{any}
valued field.

\section{$\Lambda$-metric spaces}

Throughout this section, $(\Lambda,+,\leq)$ is an ordered abelian group;
everything we need can be found in \cite{PK}.
As usual, we define $|\lambda|=\max\{\pm\lambda\}$.
\begin{Num}
We use the following slight generalization of a metric space. 
A \emph{ $\Lambda$-pseudometric} on a set $X$ is a function
$d: X\times X\too\Lambda$ which satisfies the usual axioms
of a pseudometric for all $x,y,z\in X$:
\[
d(x,x)=0\qquad
d(x,y)=d(y,x)\geq0\qquad
d(x,y)+d(y,z)\geq d(x,z).
\]
If $d(x,y)=0$ implies that $x=y$, then the $\Lambda$-pseudometric is 
called a \emph{$\Lambda$-metric}. A
$\Lambda$-pseudometric defines a topology on $X$ in the usual
way, which is Hausdorff if and only if $d$ is a $\Lambda$-metric.
For $\Lambda=(\RR,+)$, we have of course the traditional version of a
real valued (pseudo)metric.
\end{Num}
The following easy lemma will be used later.
\begin{Lem}
\label{RetrLem}
Let $d:X\times X\too\Lambda$ be a function. Suppose that for any
two elements $x,y\in X$, there exists a subset $A\subseteq X$,
containing $x$ and $y$, such that the restriction $d|_{A\times A}$
is a $\Lambda$-(pseudo)metric. If, for every such $A$,
there exists a retraction $\rho_A:X\too A$
(i.e. $\rho_A(X)=A$ and $\rho_A(a)=a$ for all $a\in A$) which diminishes
$d$ (i.e. $d(\rho_A(x),\rho_A(y))\leq d(x,y)$ for all $x,y\in X$),
then $d$ is a $\Lambda$-(pseudo)metric.

\proof
The only point to check is the triangle inequality.
So let $x,y,z\in X$, with $x,z\in A$ and $\rho_A:X\too A$ as above. Then
\[
d(x,z)=d(\rho_A(x),\rho_A(z))\leq
d(\rho_A(x),\rho_A(y))+d(\rho_A(y),\rho_A(z))
\leq d(x,y)+d(y,z).
\]
\qed
\end{Lem}
Recall that $\Lambda$ is \emph{ archimedean} if the following is true:
for any two elements $x,y>0$, there exists an $n\in\NN$ such that
$nx\geq y$. An ordered abelian group is archimedean if and only
if its only $o$-convex (\emph{order convex})
subgroups are $0$ and $\Lambda$.
Any archimedean $\Lambda$ admits an $o$-embedding into $(\RR,+)$,
which is unique up to a scaling factor;
thus, our concept of a $\Lambda$-metric is broader than the
usual concept of a metric only in the nonarchimedean case.
\begin{Num}
In the nonarchimedean case,
suppose that $\Omega\subseteq\Lambda$ is an $o$-convex subgroup, and
that $d:X\times X\too\Lambda$ is a $\Lambda$-pseudometric. Put
\[
x\approx_\Omega y\qquad\text{ if }\qquad d(x,y)\in\Omega.
\]
This is an equivalence relation on $X$, and we put
$X_x^\Omega=\{y\in X|\ x\approx_\Omega y\}$. The composite
pseudometric $X\times X\rTo^d \Lambda\too\Lambda/\Omega$ induces a
$\Lambda/\Omega$-metric on the set 
\[
\left\{X_x^\Omega|\ x\in X\right\}=X/\Omega
\]
of $\approx_\Omega$-equivalence classes in $X$. In other words,
we identify points whose distance is in $\Omega$. This
generalizes the well-known process of making a pseudometric into
a metric (the case $\Omega=0$).
\end{Num}
\begin{Num}
Given two elements $\alpha,\beta\geq0$ in $\Lambda$, we write
$\alpha\gg\beta$ if 
$\alpha>n\beta$ holds for all $n\in\NN$.
Let $\alpha\in\Lambda$ be a positive element
($\alpha>0$), and let 
\[
\Lambda^{\bra\alpha}=\{\lambda\in\Lambda|\ 
|\lambda|\leq n\alpha\text{ for some }n\in\NN\}\qquad\text{and}\qquad
\Lambda_{\bra\alpha}=\{\lambda\in\Lambda|\ 
|\lambda|\ll\alpha\}
\]
Intuitively, we truncate the $\Lambda$ behind the multiples of $\alpha$
and before $\alpha$, respectively.
Then $\Lambda_{\bra\alpha}$ is a maximal
$o$-convex subgroup of $\Lambda^{\bra\alpha}$ and the
quotient 
\[
\Lambda^{(\alpha)}=\Lambda^{\bra\alpha}/\Lambda_{\bra\alpha}
\]
is archimedean; there exists a unique $o$-homomorphism
$\phi:\Lambda^{\bra\alpha}\too\RR$ with $\phi(\alpha)=1$,
whose kernel is precisely $\Lambda_{\bra\alpha}$,
\[
0\rTo\Lambda_{\bra\alpha}\rInto\Lambda^{\bra\alpha}\rTo\RR.
\]
\end{Num}
\begin{Num}
We can make a very similar construction with any $\Lambda$-pseudometric
space $X$.
Choose a basepoint $o\in X$ and let
\[
X_o^{\bra\alpha}=\{x\in X|\ x\approx_{\Lambda^{\bra\alpha}}o\}
=\{x\in X|\ d(x,o)\leq n\alpha\text{ for some }n\in\NN\}.
\]
The composite
$\phi\circ d:X_o^{\bra\alpha}\times X_o^{\bra\alpha}\too\RR$ is
an $\RR$-pseudometric, and induces an $\RR$-metric on the
quotient 
\[
X_o^{(\alpha)}=X_o^{\bra\alpha}/\Lambda_{\bra\alpha}.
\]
\end{Num}
While all these metric concepts are rather simple and basic,
we will see that the construction of Gromov's 
\emph{asymptotic cones} \cite{Gromov} is just a special
instance of this basic method, applied to ultrapowers of metric spaces.

\section{Valuations on real closed fields}

The basic reference for this section is again \cite{PK}; the ideas
date back to Baer and Artin-Schreier.
\begin{Num}
Let $R$ be a \emph{real closed field}, i.e. an ordered field where ever
positive element is a square, and where every odd polynomial has
a zero. Let $O\subseteq R$ be an $o$-convex subring (containing
$0$ and $1$). Then $O$ is a \emph{valuation ring}: if $\alpha\in R\setminus O$,
then $\alpha^{-1}\in O$. The ring $O$ thus has a unique maximal ideal 
$M$ consisting of all nonunits. The quotient $O/M$ is again a real
closed field. If $K$ is any maximal subfield of $O$, then
$K$ is real closed and $o$-projects onto $O/M$, so
$O$ splits as a direct
sum $R=K\oplus M$ (note however that there is in general no
canonical choice for $K$),
\[
0\rTo M\rTo O\rTo^\oot O/M\rTo 0.
\]
Of course, an archimedean real closed field $R$ (such as $\RR$) contains
no proper $o$-convex subring, but there is an abundance of nonarchimedean
examples; in the next section, we will have a construction using
ultraproducts.
\end{Num}
\begin{Num}
\label{AlphaArch}
Let $R$ be a nonarchimedean real closed field, let $\alpha\gg1$, and let
$O^{\bra\alpha}=\{r\in R|\ |r|\leq\alpha^n\text{ for some }n\in\NN\}$
denote the $o$-convex subring generated by $\alpha$.
The maximal ideal is
$M^{\bra\alpha}=\{r\in R|\ |r|\leq\alpha^{-n}\text{ for all }n\in\NN\}$,
and we put 
\[
R^{\bra\alpha}=O^{\bra\alpha}/M^{\bra\alpha}.
\]
Intuitively, we have truncated the field
$R$ behind the powers of $\alpha$.
If $R=R^{\bra\alpha}$, we call $R$ an \emph{$\alpha$-archimedean field}
(because $R$ is archimedean over the subfield $\QQ(\alpha)$).
The real closed field $R^{\bra\alpha}$ is $\bar\alpha$-archimedean, where
$\bar\alpha$ is the image of $\alpha\in O^{\bra\alpha}$ in
$R^{\bra\alpha}$.
\end{Num}
\begin{Num}
\label{Log}
In general, a real closed field need not have a logarithm or exponential
function (i.e. an $o$-isomorphism $(R,+)\cong(R_{>0},\cdot)$).
In fact, no $\alpha$-archimedean real closed
field can have an exponential function.
However, we may define a \emph{formal logarithm} as follows: we take for
$(\Lambda,+)$ an $o$-isomorphic copy of the ordered abelian (multiplicative)
group $(R_{>0},\cdot)$ of positive elements. For this, we rewrite
this group additively, putting
\[
\lg:(R_{>0},\cdot)\rTo^\cong(\Lambda,+).
\]
Note that then $|\lg r|$ corresponds to $\max\{r,r^{-1}\}$.
If $O\subseteq R$ is an $o$-convex subring, we may consider the
$o$-convex subgroup 
\[
\Omega=\{\lg r|\ r>0\text{ is a unit of }O\}\subseteq\Lambda.
\]
If we put $\nu(0)=\infty$, then the map 
\[
\nu:r\mapstoo-\lg|r|+\Omega\in\Lambda/\Omega
\]
is precisely the \emph{valuation} determined by $O\subseteq R$; the
\emph{value group} is $\Gamma=\Lambda/\Omega$.
\end{Num}
\begin{Num}
Suppose that $R=R^{\bra\alpha}$ is $\alpha$-archimedean, and consider the
$o$-convex subring
\[
O_{\bra\alpha}=\{r\in R|\ |r|^n\leq\alpha\text{ for all }n\in\NN\}.
\]
Then the value group $\Gamma=\Lambda/\Omega$
is archimedean, because $\Lambda=\Lambda^{\bra{\lg\alpha}}$ and
$\Omega=\Lambda_{\bra{\lg\alpha}}$, so we have a
(nondiscrete, rank $1$) real-valued $o$-valuation $\nu$ on the field $R$.
A special instance of this construction
is Robinson's \emph{asymptotic field} $^\rho\RR$ \cite{Robinson}
\cite{Lux},
which we introduce in the next section.
\end{Num}

\section{Ultraproducts}

A basic reference for this section is \cite{CK}.
\begin{Num}
Let $I$ be an (infinite) set, and let $\mu:2^I\too\{0,1\}$ 
be a finitely additive probability measure. The collection of
all sets with $\mu$-measure $1$ is
a \emph{nonprincipal ultrafilter}. Using the axiom
of choice, the existence of such a measure can be proved
without difficulty. (Note that every subset of $I$ is required to be
$\mu$-measurable! Clearly, every finite set has measure $0$,
and every cofinite set has measure $1$.)

Suppose that $(\mathcal Q_i)_{i\in I}$ is a family of
first order structures: a family of groups, rings, fields,
or metric spaces. The direct product $\prod\mathcal Q_i$
of these structures will in general be of a weaker type;
for example, the direct product of fields is a ring
with zero-divisors. This deficiency can be corrected
using the measure $\mu$: two elements
$(q_i)_{i\in I},(q'_i)_{i\in I}$ of the direct product
are identified if their difference set $\{j\in I|\ q_j\neq q_j'\}$
has $\mu$-measure $0$. The resulting collection of equivalence classes
is the \emph{ultraproduct}
\[
\prod_\mu\mathcal Q_i.
\]
An ultraproduct has the same kind of first-order properties as
its factors; for example, an ultraproduct of fields is again a field.
Note that this construction is completely analogous
to the construction of the Hilbert space $L^2$, where integrable
functions are identified if they differ on a set of Lebesgue
measure $0$.
\end{Num}
\begin{Num}
An important special case is when all structures $\mathcal Q_i$
are equal to one fixed structure $\mathcal Q$; in this case, one obtains
an \emph{ultrapower} 
\[
^*\mathcal Q=\prod_\mu\mathcal Q
\]
which is also called a \emph{nonstandard model} of $\mathcal Q$.
If $\mathcal Q=(\RR,+,\cdot,\leq)$
is the field of real numbers, the resulting field
$\prod_\mu(\RR,+,\cdot,\leq)=(^*\RR,+,\cdot,\leq)$ is the field of
\emph{nonstandard real numbers}.
This is a nonarchimedean real closed field.
\end{Num}
\begin{Num}
Let $\alpha\in{}^*\RR$ be an infinitely large nonstandard real,
and let $\rho=\alpha^{-1}$.
The fields
\[
^\rho\RR=(^*\RR)^{\bra\alpha}
\]
were first studied by Robinson \cite{Robinson} \cite{Lux},
and are sometimes called
\emph{Robinson's asymptotic fields}.
\end{Num}
\begin{Num}
Suppose now that each $\mathcal Q_i=(d_i:X_i\times X_i\too\RR)$ is
an $\RR$-metric space. Then the ultraproduct
\[
\prod_\mu\mathcal (d_i:X_i\times X_i\too\RR)=
(\mathbf d:\mathbf X\times\mathbf X\too{}^*\RR)
\]
is a $\Lambda$-metric space, with $\Lambda=(^*\RR,+)$.
If we pick a point $o\in\mathbf X$, a number $\alpha\gg1$ in
${}^*\RR$, then 
\[
X^{\bra\alpha}_o=\left\{x\in\mathbf X\left|\ \alpha^{-1}\mathbf d(x,o)
\leq n
\text{ for some }n\in\NN\right\}\right.
\]
and the $\RR$-metric space
\[
\mathbf X^{(\alpha)}_o=\mathbf X^{\bra\alpha}_o/\Lambda_{\bra\alpha}
\]
is precisely Gromov's
\emph{asymptotic cone} of the family $(X_i)_{i\in I}$, where all points
$x,y$ with infinitesimal distance $\alpha^{-1}\mathbf d(x,y)$ are
identified \cite{Gromov}.
The general construction is due to Van den Dries and Wilkie
\cite{VDW}.
\end{Num}

\section{Riemannian symmetric spaces over real closed fields}

For the geometry of Riemannian symmetric spaces, see \cite{Eberlein},
\cite{Helgason}, \cite{BGS}, and \cite{BH}.
\begin{Num}
Let $P_n$ denote the set of all symmetric positive
definite real $n\times n$-matrices with determinant $1$. Every such matrix
is of the form $X=gg^T$, for some $g\in\SL_n\RR$. The set $P_n$ can
be identified with the Riemannian symmetric space $SL_n\RR/\SO(n)$
consisting of all elliptic polarities of real projective $n-1$-space.
The group $\SL_n\RR$ acts as 
\[
(g,X)\mapstoo gXg^T
\]
(the geodesic reflection at $X\in P_n$ is $Y\mapstoo XY^{-1}X$).
Up to a real scaling factor, the metric $d_R$ induced by the
(unique invariant) Riemannian metric of $P_n$ is given as
follows. For $X\in P_n$,
let $Y=\log(X)$ denote the unique traceless symmetric matrix
with $\exp(Y)=X$. Then 
\[
d_R(\mathbf 1,X)^2=\tr(Y^2).
\]
\end{Num}
\begin{Num}
\label{MetricOnX}
Now $P_n$ is an algebraic variety defined over $\QQ$, and we may 
consider its set of $R$-points $P_n(R)$
over any real closed field $R$.
The metric $d_R$, however, involves an analytic function, the
logarithm, which need not be defined over an arbitrary real closed field
$R$.
We fix this problem, using the formal logarithm 
$\lg:R_{>0}\too\Lambda$ defined in \ref{Log}.
Given a diagonal matrix $X=\mathrm{diag}(x_1,\ldots,x_n)\in P_n(R)$,
we define the \emph{$\Lambda$-valued distance}
\[
d(\mathbf 1,X)=|\lg x_1|+\cdots+|\lg x_n|\in\Lambda.
\]
This function is obviously invariant under coordinate permutations
(the Weyl group action for $\SL_nR$); it follows that we can use the
$\SL_n R$-action to
extend $d$ to a well-defined distance function on $P_n(R)$, setting
\[
d(gg^T,gXg^T)=d(\mathbf 1,X)
\]
(where $X$ is a diagonal matrix).
On the set $A$ of all diagonal matrices in $P_n(R)$, this is clearly a
$\Lambda$-metric (essentially, it is the Manhattan Taxi Metric).
By Lemma \ref{RetrLem} and Kostant's Convexity Theorem \cite{Kostant},
$d$ is indeed an $\SL_nR$-invariant metric on $P_n(R)$.
Over the reals, we may take $\Lambda=\RR$ and $\lg=\log$; then
$d$ and $d_R$ are different, but (bi-Lipschitz) equivalent metrics.
\end{Num}
\begin{Num}
We have successfully made $X=P_n(R)$ into a $\Lambda$-metric space.
Suppose now that $O\subsetneq R$ is an $o$-convex subring.
Put as before
$\Omega=\{\lg r| r>0\text{ and }r\text{ is a unit of }O\}\subseteq\Lambda$.
Then $\Gamma=\Lambda/\Omega$ is the value group of the
valuation determined by $O$.
The group $\SL_nR$ acts by isometries on the $\Lambda/\Omega$-metric
space $P_n(R)/\Omega$, and it is not difficult to check the following
result:
\[
P_n(R)/\Omega=\SL_n(R)/\SL_n(O).
\]
\end{Num}
\begin{Num}
Suppose now that $R=R^{\bra\alpha}$ is $\alpha$-archimedean,
see \ref{AlphaArch}, and that $O=O_{\bra\alpha}$. Then the value group
$\Gamma=\Lambda/\Omega=\Lambda^{(\alpha)}$ is archimedean
(we could take for example $R={}^\rho\RR$,
Robinson's asymptotic field).
Consequently,
\[
P_n(R)/\Omega=\SL_n(R)/\SL_n(O)
\]
is an $\RR$-metric space.
The main result of Section \ref{LambdaBuildingsSection}
is that this quotient is a
(nondiscrete) affine $\RR$-building
(of type $\widetilde A_{n-1}$, in fact).
The group $\SL_n(R)$ acts as a
transitive automorphism group on its vertices.
More generally, if $O\subseteq R$ is any $o$-convex valuation ring,
then $\SL_n(R)/\SL_n(O)$ is an affine $\Lambda/\Omega$-building
in the sense of Bennett \cite{Bennett}.
\end{Num}
\begin{Num}
In this section, we have concentrated on the Riemannian symmetric space
$\SL_n\RR/\SO(n)$. Everything we have done works in the same generality
for any  Riemannian symmetric space of noncompact type. In fact,
any irreducible Riemannian symmetric space can be equivariantly
embedded in $P_n$ as a totally geodesic submanifold;
the embedding can be chosen to be algebraic over $\QQ$.
Then one obtains without much difficulty the following result.
\end{Num}
\begin{Thm}
Let $R$ be a (nonarchimedean) real closed field, let
$O$ be an $o$-convex subring,
let $X$ be a Riemannian symmetric space of noncompact type. Let
$G$ be the corresponding semisimple real Lie group,
with maximal compact subgroup $K$. Then we can view $G$ and $K$ as
real algebraic groups defined over $\QQ$, and $X=G/K$ as a real algebraic
variety. Define $\Lambda$ and $\Omega$ as above.
Then there exists a $G(R)$-invariant $\Lambda$-metric on $X(R)=G(R)/K(R)$
(semialgebraic over $\QQ$), and there
is a natural equivariant identification
\[
X(R)/\Omega=G(R)/G(O)
\]
of $\Lambda/\Omega$-metric spaces.
If $R$ is in addition $\alpha$-archimedean, and $O=O_{\bra\alpha}$,
then $X(R)/\Omega$ is an $\RR$-metric space.
\qed
\end{Thm}
Now we consider as a special case the ultrapower
$\prod_\mu(X\times X\rTo\RR)=(\mathbf X\times\mathbf X\too{}^*\RR)$.
Pick a basepoint $o\in\mathbf X$ and an infinitely large nonstandard
real $\alpha\gg 1$. Let $\rho=1/\alpha$.
\begin{Cor}
\label{HomogeneousCone}
For the asymptotic cone of the Riemannian symmetric space $X$, we have
\[
\mathbf X^{(\alpha)}_o=G(R)/G(O),
\]
where $R={}^\rho\RR=(^*\RR)^{\bra\alpha}$ is Robinson's asymptotic field,
and $O=O_{\bra{\bar\alpha}}$.
\qed
\end{Cor}

\section{Affine $\Lambda$-buildings}
\label{LambdaBuildingsSection}
\begin{Num}
Suppose that $W$ is a finite Coxeter group satisfying the crystallographic
condition. Then $W$ acts naturally on a $\ZZ$-lattice $L\cong\ZZ^n$
in $\RR^n$.
Let $\Lambda$ be an ordered abelian group; then $W$ acts
naturally on $A=L\otimes\Lambda\cong\Lambda^n$. Let
$\overline W=W\ltimes(A,+)$ denote 
the corresponding affine reflection group.
Using the $W$-invariant inner product $L\otimes L\too\ZZ$, it is
possible to construct a $\overline W$-invariant $\Lambda$-metric
on $d:A\times A\too\Lambda$. (For the case of the symmetric group,
the metric  -- the Manhattan Taxi Metric -- was given in 
\ref{MetricOnX} above.)
Using the action of $W$ on $A$, it makes sense to talk about affine
reflection hyperplanes and closed halfspaces. A subset $B\subseteq A$
is called \emph{$W$-convex} if it is the intersection of finitely many
halfspaces.
\end{Num}
\begin{Num}
Suppose that $X$ is a $\Lambda$-metric space. An
\emph{atlas} on $X$
is a collection $\mathcal A$ of $\Lambda$-isometric injections
$\phi:A\too X$, called \emph{coordinate charts},
with the following properties.

\textbf{(A1)} If $\phi$ is in $\mathcal A$ and $w\in\overline W$, then
$\phi\circ w:A\too X$ is in $\mathcal A$.

\textbf{(A2)} Given two charts $\phi_1,\phi_2$, the set
$B=\phi^{-1}(\phi(A))$ is $W$-convex, and there exists a $w\in\overline W$
with $\phi_1|_B=\phi_2\circ w|_B$.

The sets $F=\phi(A)$ are called \emph{apartments}; the image $S=\phi(S_0)$
of the basic Weyl cone $S_0\subseteq A$ is called a \emph{sector}.

\textbf{(A3)} Given $x,y\in X$, there exists an apartment
$F=\phi(A)$ containing $x$ and $y$.

\textbf{(A4)} Given two sectors $S_1,S_2\subseteq X$, there exist
subsectors $S_1'\subseteq S_1$ and $S_2'\subseteq S_2$ and an apartment
$F$ containing $S_1'\cup S_2'$.

\textbf{(A5)} If $F_1,F_2,F_3$ are apartments such that each of the three
sets $F_i\cap F_j$, $i\neq j$, is a halfapartment
(i.e. the $\phi$-image of a halfspace),
then $F_1\cap F_2\cap F_3\neq\emptyset$.

\textbf{(A6)}
For any apartment $F$ and any $x\in F$, there exists a retraction
$\rho_{x,F}:X\too F$ which diminishes distances, with
$\rho_{x,F}^{-1}(x)=\{x\}$.

\smallskip\noindent
The pair $(X,\mathcal A)$ is called an \emph{affine $\Lambda$-building}
\cite{Bennett}.
The atlas $\mathcal A$ is, in general, not unique, but it is always contained
in a unique maximal atlas.
\end{Num}
If we take for $X$ a Riemannian symmetric space and the maximal flats
as the apartments, then $X$ satisfies axioms \textbf{(A1)} and \textbf{(A3)},
and, rather trivially, axioms \textbf{(A2)} and \textbf{(A5)}, while
axiom \textbf{(A4)} is only approximately true.
\begin{Num}
For $n=1$, an affine $\Lambda$-building is the same as a $\Lambda$-tree
\cite{Bennett} \cite{AB};
affine $\Lambda$-buildings are in a sense higher-dimensional versions of 
($\Lambda$-)trees.
The notion of affine $\Lambda$-buildings is due to Bennett \cite{Bennett};
affine $\RR$-buildings are
the same as a Tits' \emph{syst\`emes d'appartements} \cite{Tits}.
Kleiner-Leeb \cite{KL} give a different set of axioms for spaces they
call \emph{euclidean buildings}; Parreau \cite{Parreau} proved that
their spaces are special cases of affine $\RR$-buildings.
Finally, there are the affine buildings
in the proper sense of building theory \cite{Brown} \cite{Ronan}.
One has the following inclusions:
{\small\[
\{\text{affine $\Lambda$-buildings}\}\supsetneq
\{\text{affine $\RR$-buildings}\}\supsetneq
\{\text{Kleiner-Leeb euclidean buildings}\}\supsetneq
\{\text{affine buildings}\}
\]}%
\end{Num}
\begin{Num}
\label{BuildingAtInfinity}
Let $X$ be an affine $\Lambda$-building, let $o\in X$. One can construct two
spherical buildings from $X$. Consider the collection $\mathbf{Sec}_o$ of all
sectors $\phi(S_0)$ whose tip $\phi(0)$ is $o$
(recall that $S_0\subseteq L\otimes\Lambda$ is the basic Weyl cone).
Then $\mathbf{Sec}_o$ generates in a rather natural way a poset, the
\emph{spherical building at infinity}. Using the axioms above, one can
show that different basepoints lead to canonically isomorphic
buildings at infinity, and we denote the resulting building by
\[
\Delta_\infty^\mathcal A X.
\]
This building depends on the atlas $\mathcal A$;
if $\mathcal A$ is maximal, we omit it and write $\Delta_\infty X$
\end{Num}
\begin{Num}
\label{LocalBuilding}
One can construct another spherical building from $\mathbf{Sec}_o$;
here, we identify two sectors
$S_1,S_2\in\mathbf{Sec}_o$ if they agree inside an open ball around $o$
(with respect to the $\Lambda$-metric).
The corresponding spherical building is denoted $\Delta_oX$; it is independent
of the apartment system (because it is defined by local data). There
is a canonical building epimorphism
\[
\Delta_\infty^\mathcal A X\too\Delta_oX.
\]
\end{Num}
\begin{Num}
Let $R$ be a real closed field and let $O\subsetneq R$ be an
$o$-convex subring. Let $X=G/K$ be a Riemannian symmetric space of
noncompact type. Our first main result is as follows.
\end{Num}
\begin{Thm}
The quotient $G(R)/G(O)$ is an affine $\Gamma$-building, where
$\Gamma=\Lambda/\Omega$ is the value group of the valuation determined
by $O$, see \ref{Log}. In particular, if
$R=R^{\bra\alpha}$ is $\alpha$-archimedean and $O=O_{\bra\alpha}$,
then $G(R)/G(O)$ is
an affine $\RR$-building. The automorphism group of the building at
infinity (with respect to the apartment system we construct) is
the group $G(R)$, provided that no simple factor of $G$ has
rank $1$.
\end{Thm}
The proof depends on various classical results about semisimple Lie
groups. As the apartment system, we take the images of the maximal
flats in the Riemannian symmetric space $G(R)/K(R)$ under the canonical map
\[
G(R)/K(R)=G(R)/K(O)\too G(R)/G(O).
\]
Let $G=KAU$ be an Iwasawa decomposition.
For example, axiom \textbf{(A3)} is a consequence of the $KAK$-decomposition,
$G(R)=K(R)A(R)K(R)$, while \textbf{(A2)} follows from Kostant's Convexity
Theorem. The proof of \textbf{(A6)} depends also on the
Convexity Theorem, and the Iwasawa projection $g=kau\mapstoo a$.
\begin{Cor}
In the special case of Robinson's asymptotic fields $^\rho\RR$, we obtain
in this way the full (maximal) apartment system. Thus, if $G$ has no 
simple factor of rank $1$, then $G(^\rho\RR)$ is the full automorphism
group of $X$.
\end{Cor}
The corollary uses the model-theoretic
fact that ultrapowers are saturated.

\section{Topological rigidity of affine $\RR$-buildings}

Let $\Delta$ be the underlying simplicial complex of a (thick)
rank $k$ building \cite{Brown} \cite{Ronan}, and
let $|\Delta|$ denote its geometric realization (topologized in
any sensible way). It is not difficult to recover the combinatorial
structure of $\Delta$ from the topological space $|\Delta|$;
the $k-1$-skeleton of $|\Delta|$ is precisely the set of all points
in $|\Delta|$ which do not have a locally euclidean neighborhood.
It follows that buildings whose geometric realizations are homeomorphic
are isomorphic.
This is considerably more involved for nondiscrete affine $\RR$-buildings,
since here, it may happen that no point has a locally euclidean neighborhood.
Topological rigidity of affine $\RR$-buildings was proved first by
Kleiner-Leeb \cite{KL}.
The appearance of local homology groups in their argument
suggests that there should be a simple sheaf-theoretic proof.
This is indeed the case. The amount of sheaf theory needed here is
very modest. The first introductory pages of \cite{Bredon} completely
suffice; in particular, we do \emph{not} need such sophisticated
tools as sheaf-theoretic (co)homology, although some ideas are
inspired by \cite{Lowen}.
\begin{Num}
To any Hausdorff space $X$, one can associate two basic presheaves:
firstly, the graded $\ZZ/2$-module valued
presheaf $U\mapstoo H_\bullet(X,X\setminus U;\ZZ/2)$
(ordinary singular homology with $\ZZ/2$-coefficients)
and secondly, the presheaf $U\mapstoo 2^U$ which assigns to $U$
the boolean algebra of all subsets of $U$. Let
\[
\mathcal H_\bullet=\mathit{Sheaf}(U\mapstoo H_\bullet(X,X\setminus U;\ZZ/2)
\quad\text{ and }\quad
\mathcal S=\mathit{Sheaf}(U\mapstoo 2^U)
\]
denote the corresponding sheaves.
The stalks of the first sheaf are the local homology groups,
\[
\mathcal H_{\bullet,x}=H_\bullet(X,X\setminus\{x\};\ZZ/2)
\]
(this follows from the axiom of compact supports for singular homology),
while the stalk $\mathcal S_x$ consists of germs of subsets of $X$ near $x$.
Note that the closure operation for subsets is compatible with
restriction and thus induces a closure operation on
germs of subsets near $x$. Given a local section
$s:U\too\mathcal H_\bullet(U)$
of $\mathcal H_\bullet$ over $U$, there is the closed subset
$\supp(s)=\{x\in U|\ s_x\neq 0\}$.
This yields a natural map $ssg:\mathcal H_{\bullet,x}\rTo\mathcal S_x$
on the stalks: any element $\xi\in H_\bullet(X,X\setminus\{x\};\ZZ/2)$
determines a subset germ $ssg(\xi)$ at $x$.
For $\xi,\eta\in\mathcal H_{\bullet,x}$, we define
\[
\xi*\eta=ssg(\xi)\cap ssg(\eta)\setminus ssg(\xi+\eta)\in\mathcal S_x.
\]
\end{Num}
\begin{Num}
Suppose now that $X$ is an affine $\RR$-building. Let $o\in X$ and
recall from
\ref{BuildingAtInfinity} and
\ref{LocalBuilding}
the definition of $\mathbf{Sec}_o$ and the building $\Delta_oX$.
Each sector $S\in\mathbf{Sec}_o$ determines a germ $l(S)\in\mathcal{S}_o$;
the resulting sub-poset of $\mathcal S_o$ generated by the intersections
of these germs is precisely $\Delta_o X$,
\[
\Delta_oX\rInto\mathcal S_o.
\]
Now let $F$ be an apartment of $X$ containing $o$. Using the retraction
$\rho_{o,F}:X\too F$, one sees that $F$ determines an element
$\xi_F\in H_\bullet(X,X\setminus\{o\};\ZZ/2)$, given by the image of
$H_\bullet(F,F\setminus\{o\};\ZZ/2)$.
To simplify the argument, we assume that all the buildings $\Delta_oX$
are \emph{thick}; in the setting of the
proof of the Margulis conjecture, this is the case.
(The general case, where the $\Delta_oX$ may be weak buildings requires
some minor refinements.)
\end{Num}
\begin{Prop}
There is an isomorphism
\[
H_{\bullet+1}(X,X\setminus\{o\};\ZZ/2)\cong
H_\bullet(X\setminus\{o\};\ZZ/2)\cong H_\bullet(|\Delta_oX|;\ZZ/2).
\]

\noindent\em
The first isomorphism comes from the contractibility of $X$.
The main point of the proof of the second isomorphism
is to show that every element
$\xi\in H_\bullet(X\setminus\{o\};\ZZ/2)$ arises as above from
(finitely many)
apartments. Here, the key idea is to use the fact that every chamber
in $\Delta_oX$ has a small neighborhood which lifts isometrically
into the neighborhood of some chamber of $\Delta_\infty X$.
\qed
\end{Prop}
Using this result, we can characterize the germs
$l(S)\in \mathcal S_o$.
\begin{Thm}
Every chamber $l(S)$ of $\Delta_oX$ is of the form $\overline{\xi*\eta}$,
for suitably chosen $\xi,\eta$. Every element $\xi*\eta$
which is minimal
(i.e. there exist no $\xi',\eta'$ with
$\emptyset\neq \xi'*\eta'\subsetneq\xi*\eta$)
represents a chamber $l(S)$ of $\Delta_oX$.
\qed
\end{Thm}
\begin{Cor}
Let $f:X\rTo X'$ be a homeomorphism of affine $\RR$-buildings.
Then $f$ induces an isomorphism between $\Delta_oX$ and
$\Delta_{f(o)}X'$, for every $o\in X$.
\end{Cor}
From this local result, one can without difficulty derive the following.
\begin{Thm}
Let $s:X\rTo \mathcal H_\bullet$ be a global section of the sheaf
$\mathcal H_\bullet\rTo X$ over the affine $\RR$-building $X$.
Then $A=\supp(s)\subset X$ is an apartment if and only if
(1) $A$ is homeomorphic to $\RR^n$, for some $n$, and
(2) $ssg(s_a)$ is an apartment in $\Delta_a$ for every $a\in A$.
\end{Thm}
\begin{Cor}
Let $f:X\rTo X'$ be a homeomorphism of affine $\RR$-buildings. Then
$f$ carries apartments to apartments.
\end{Cor}
\begin{Cor}
\label{TopologicalRigidity}
A homeomorphism of affine $\RR$-buildings $X,X'$ induces an
isomorphism (in a sense to be made precise); moreover,
$\Delta_\infty X\cong\Delta_\infty X'$.
\end{Cor}

\section{Proof of the Margulis conjecture}

Now we outline a proof of the Margulis conjecture, using the previous
results.
\begin{Num}
Let $X,X'$ be Riemannian symmetric spaces of noncompact type,
without de Rham factor of rank $1$, and let 
\[
f:X\rTo X
\]
be an $(L,C)$-quasi-isometry,
\[
L^{-1}\cdot d(x,y)-C\leq d'(f(x),f(y))\leq L\cdot d(x,y)+C.
\]
Taking ultrapowers of $X$ and $X'$, we obtain
an $(L,C)$ quasi-isometry
\[
^*f:{}^*X\too{}^*X'
\]
between the $\Lambda$-metric spaces $^*X,{}^*X'$, where
$(\Lambda,+)\cong(^*\RR_{>0},\cdot)$ (we use the metric $d$ introduced
in \ref{MetricOnX},
which is bi-Lipschitz equivalent to the metric $d_R$ induced
by the Riemannian metric).
\end{Num}
\begin{Num}
Pick any point $o\in {}^*X$, and an element $\alpha\in{}^*\RR$, with
$\alpha\gg 1$. The properties of a quasi-isometry and the finiteness of
the constants $L,C$ imply that we have an $(L,C)$-quasi-isometry
\[
(^*X)_o^{\bra\alpha}\rTo^{^*f} (^*X')_{f(o)}^{\bra\alpha}
\]
of $\Lambda^{\bra\alpha}$-metric spaces,
which descends to a (bi-Lipschitz) homeomorphism
\[
(^*X)_o^{(\alpha)}\rTo^F (^*X')_{f(o)}^{(\alpha)}
\]
on the asymptotic cones.
\end{Num}
\begin{Num}
Let $R={}^\rho\RR=(^*\RR)^{\bra\alpha}$
denote Robinson's asymptotic field, and
let $O=O_{\bra{\bar\alpha}}$. By the results in Section 
\ref{LambdaBuildingsSection},
$Z=(^*X)_o^{(\alpha)}=G(R)/G(O)$ and
$Z'=G'(R)/G'(O)=(^*X')_{f(o)}^{(\alpha)}$ are affine $\RR$-buildings.
By \ref{TopologicalRigidity}, we have an isomorphism 
\[
\Delta_\infty Z\cong\Delta_\infty Z',
\]
whence a group isomorphism $G(^\rho\RR)\cong G'(^\rho\RR)$.
This implies that there is a Lie group isomorphism $G\cong G'$,
and so, an equivariant diffeomorphism
\[
X=G/K\cong G'/K'=X'.
\]
It follows that the metrics on the de Rham factors of $X$ can be rescaled
such that $X$ and $X'$ are isometric.
\end{Num}
\begin{Num}
A careful analysis of the quasi-isometry shows that a stronger result holds,
see \cite{KL}:
there exists a (necessarily unique) isometry $\bar f:X\too X'$ such that
$\bar f$ has bounded distance from $f$.
For the proof, one has to consider different choices for
$\alpha$ and the base point $o$; the main step is to show that
$f$ is (after rescaling the metrics on the de Rham factors) a coarse
isometry.
The original proof of this by Kleiner-Leeb can be simplified in a substantial
way, using $\Lambda$-buildings and model theoretic methods.
\end{Num}


\begin{thebibliography}{XX}

\bibitem{AB}
R. Alperin and H. Bass,
Length functions of group actions on $\Lambda$-trees. 
\emph{Combinatorial group theory and topology} (Alta, Utah, 1984), 
265--378,
Princeton Univ. Press, Princeton, NJ, 1987.

\bibitem{BGS}
W. Ballmann, M. Gromov and V. Schroeder,
\emph{Manifolds of nonpositive curvature}.
Birkh\"auser Boston, Inc., Boston, MA, 1985.

\bibitem{Bennett}
C. Bennett,
Affine $\Lambda$-buildings. I.
Proc. London Math. Soc. 68 (1994) 541--576.

\bibitem{Bredon}
G. Bredon,
\emph{Sheaf theory}, 2nd ed.
Springer-Verlag, New York, 1997.

\bibitem{BH}
M. Bridson and A. Haefliger,
\emph{Metric spaces of non-positive curvature}.
Springer-Verlag, Berlin, 1999.

\bibitem{Brown}
K. Brown,
\emph{Buildings.}
Springer-Verlag, New York (1989).

\bibitem{CK}
C.C. Chang and H.J. Keisler,
\emph{Model theory}. 3rd ed.
North-Holland Publishing Co., Amsterdam, 1990.

\bibitem{Eberlein}
P. Eberlein,
\emph{Geometry of nonpositively curved manifolds}.
University of Chicago Press,
Chicago, IL, 1996.

\bibitem{EF}
A. Eskin and B. Farb,
Quasi-flats and rigidity in higher rank symmetric spaces.
J. Amer. Math. Soc. 10 (1997) 653--692.

\bibitem{Gromov}
M. Gromov,
Asymptotic invariants of infinite groups.
\emph{Geometric group theory, Vol. 2 (Sussex, 1991)}, 1--295,
Cambridge Univ. Press, Cambridge, 1993. 

\bibitem{Helgason}
S. Helgason,
\emph{Differential geometry, Lie groups, and symmetric spaces}.
Academic Press, Inc., New York-London, 1978.

\bibitem{KL}
B. Kleiner and B. Leeb,
\emph{Rigidity of quasi-isometries for symmetric spaces and
Euclidean buildings.}
Inst. Hautes \'Etudes Sci. Publ. Math. No. 86 (1997) 115--197. 

\bibitem{Kostant}
B. Kostant,
On convexity, the Weyl group and the Iwasawa decomposition. 
Ann. Sci. \'Ecole Norm. Sup. (4) 6 (1973), 413--455. 

\bibitem{Lowen}
R. L\"owen,
Topology and dimension of stable planes: on a conjecture of H. Freudenthal.
J. Reine Angew. Math. 343 (1983) 108--122.

\bibitem{Lux}
W. Luxemburg,
On a class of valuation fields introduced by A. Robinson.
Israel J. Math. 25 (1976) 189--201.

\bibitem{Mostow}
G.D. Mostow,
\emph{Strong rigidity of locally symmetric spaces.}
Princeton University Press,
Princeton, N.J., 1973.

\bibitem{Parreau}
A. Parreau,
Immeubles affines: construction par les normes et \'etude des isom\'etries.
\emph{Crystallographic groups and their generalizations (Kortrijk, 1999)},
263--302, Contemp. Math., 262, Amer. Math. Soc., Providence, RI, 2000.

\bibitem{PK}
S. Prie\ss-Crampe,
\emph{Angeordnete Strukturen: Gruppen, K\"orper, projektive Ebenen.}
Springer-Verlag, Berlin, 1983.

\bibitem{Robinson}
A.H. Lightstone and A. Robinson,
\emph{Nonarchimedean fields and asymptotic expansions}.
North-Holland Publishing Co.,
New York, 1975.

\bibitem{Ronan}
M. Ronan,
\emph{Lectures on buildings}.
Academic Press, Inc., Boston, MA, 1989.

\bibitem{Thornton}
B. Thornton,
\emph{Asymptotic cones of symmetric spaces},
PhD Thesis, Univ. Utah, 2002.

\bibitem{Tits}
J. Tits,
Immeubles de type affine.
\emph{Buildings and the geometry of diagrams (Como, 1984)}, 159--190, 
Lecture Notes in Math. 1181, 
Springer, Berlin, 1986. 

\bibitem{VDW}
L. van den Dries and A. Wilkie, 
Gromov's theorem on groups of polynomial growth and elementary logic.
J. Algebra 89 (1984) 349--374.

\end{thebibliography}
\end{document}